\documentclass[reqno,12pt]{amsart}

\usepackage{latexsym}
\usepackage[centertags]{amsmath}
\usepackage{amsfonts}
\usepackage{amssymb}
\usepackage{amsthm}
\usepackage{newlfont}
\usepackage[square]{natbib}

\def\mytopsep{3mm}

\newtheoremstyle{myplain}{\mytopsep}{\mytopsep}{\itshape}{0pt}{\bfseries}{.}{3mm}{}
\newtheoremstyle{mydefinition}{\mytopsep}{\mytopsep}{\normalfont}{0pt}{\bfseries}{.}{3mm}{}

\newtheoremstyle{myremark}{\mytopsep}{\mytopsep}{\normalfont}{0pt}{\bfseries}{.}{3mm}{}

\theoremstyle{myplain}

\newtheorem{thm}{Theorem}[section]
\newtheorem{cor}[thm]{Corollary}
\newtheorem{lem}[thm]{Lemma}
\newtheorem{prop}[thm]{Proposition}

\theoremstyle{mydefinition}
\newtheorem{dfn}[thm]{Definition}
\theoremstyle{myremark}
\newtheorem{rem}[thm]{Remark}
\newtheorem{exa}[thm]{Example}

\allowdisplaybreaks[1]

\makeatletter\@addtoreset{equation}{section}\makeatother

\def\bcup{\bigcup\limits}
\def\eql{\simeq}

\def\supp{\mathop{\mbox{supp}}}
\def\ord{\mathrm{ord}}
\def\mb{\mathbf}
\def\diag{\mathrm{diag }}

\def\fra{{\mathrm{fra}}}

\def\QQ{\mathbb{Q}}

\def\xx{x^{-1}}

\def\CC{\mathbb{C}}
\def\ZZ{\mathbb{Z}}
\def\RR{\mathbb{R}}

\def\res{\mathop{\mathrm{Res}}}
\def\ct{\mathop{\mathrm{CT}}}

\def\CT{\mathop{\mathrm{CT}}}

\newcommand{\pad}[2]{\displaystyle\frac{\partial #2}{\partial #1}}

\renewcommand{\ll}{\langle\!\langle}
\renewcommand{\gg}{\rangle\!\rangle}

\begin{document}

\title{A Residue Theorem for Malcev-Neumann Series}

\author{Guoce Xin}
\address{Department
of Mathematics\\
Brandeis University\\
Waltham MA 02454-9110} \email{guoce.xin@gmail.com}

\date{February 21, 2005}
\begin{abstract}
In this paper, we establish a residue theorem for Malcev-Neumann
series that requires few constraints, and includes previously
known combinatorial residue theorems as special cases. Our residue
theorem identifies the residues of two formal series that are
related by a change of variables. We obtain simple conditions for
when a change of variables is possible, and find that the two
related formal series in fact belong to two different fields of
Malcev-Neumann series. The multivariate Lagrange inversion formula
is easily derived and Dyson's conjecture is given a new proof and
generalized.

\end{abstract}
\maketitle


\section{Introduction}

Let $K$ be a field. Jacobi \cite{jac} used the ring $K((x_1,\dots
,x_n))$ of Laurent series, formal series of monomials where the
exponents of the variables are bounded from below, to give the
following residue formula.

\begin{thm}[Jacobi's Residue Formula]
Let $f_1(x_1,\dots ,x_n), \ldots , f_n(x_1,\dots,x_n)$ be Laurent
series. Let $b_{i,j}$ be integers such that
$f_i(x_1,\dots ,x_n)/x_1^{b_{i,1}}\cdots x_n^{b_{i,n}}$ is a
formal power series with nonzero constant term. Then for any
Laurent series $\Phi(y_1,\dots ,y_n)$, we have
\begin{equation}
\res_{x_1,\dots ,x_n} \left| \pad{x_{j}} {f_i}\right|_{1\le i,j\le
n} \Phi(f_1,\dots ,f_n) =\left| b_{i,j}\right|_{1\le i,j\le n}
\res_{y_1,\dots ,y_n} \Phi(y_1,\dots ,y_n),
\end{equation}
where $\res_{x_1,\dots ,x_n}  $ means to take the coefficient of
$x_1^{-1}\cdots x_n^{-1}$.
\end{thm}
Note that the convergence of $\Phi(f_1,\dots,f_n)$ is obviously
required.

Jacobi's residue formula is a well-known result in combinatorics.
It equates the residues of two formal series related by a change
of variables. It has many applications and has been studied by
several authors, e.g., Goulden and Jackson
\cite[p.~19--22]{g-jackson}, and Henrici \cite{henrici}. However, Jacobi's
formula is rather restricted in application for two reasons: the
conditions on the $f_i$ are too strong, and the condition on
$\Phi$ is not easy to check: given $f_i$, when does
$\Phi(f_1,\dots ,f_n)$ converge?

We can obtain different residue formulas by considering different
rings containing the ring of formal power series
$K[[x_1,\dots,x_n]]$. In obtaining such a formula, we usually
embed $K[[x_1,\dots,x_n]]$ into a ring or a field consisting of
formal Laurent series, but the embedding is not unique in the
multivariate case. Besides Jacobi's residue formula, Cheng et al.\
\cite{reversion} studied the ring $K_h((x_1,\dots ,x_n))$ of
homogeneous Laurent series (formal series of monomials whose total
degree is bounded from below), and used homogeneous expansion to
give a residue formula. But the above restrictions still exist for
the same reason. We will use a more general setting to avoid the
above problems.

Let $\mathcal{G}$ be a \emph{totally ordered group}, i.e., a group
with a total ordering $\le$ that is compatible with its group
structure. Let $K_w[\mathcal{G}]$ be the set of
\emph{Malcev-Neumann series} (MN-series for short) on
$\mathcal{G}$ over $K$ relative to $\le$: an element in
$K_w[\mathcal{G}]$ is a series $\eta=\sum_{g\in \mathcal{G}}a_gg$
with $a_g\in K$, such that the support $\{\, g\in \mathcal{G}:
a_g\ne 0\, \}$ of $\eta$ is a well-ordered subset of
$\mathcal{G}$.

By a theorem of Malcev \cite{malcev}
and Neumann
\cite{neumann} (see also \cite[Theorem~13.2.11]{passmann}),
$K_w[\mathcal{G}]$ is a division algebra that includes the group algebra
$K[\mathcal{G}]$ as a subalgebra. We study the field of MN-series
on a totally ordered abelian group, and show that the field of
iterated Laurent series $K\ll x_1,\dots ,x_n\gg$, which has been
studied in \cite[Chapter~2]{xinthesis}, is a special kind of
MN-series.

We obtain a residue theorem for $K_w[\mathcal{G}\oplus \ZZ^n]$,
where $x_1,\dots ,x_n$ represent the generators of $\ZZ^n$. This
new residue formula includes the previous residue theorems of
Jacobi and Cheng et al.\ as special cases. It is easier to apply
and more general: the conditions on the $f_i$ are dropped since we
are working in a field; the condition on $\Phi$ is replaced with
a simpler one and we find that the two related formal series in
fact belong to two different fields of MN-series. In particular,
our theorem applies to any rational function $\Phi$.

In section 2 we review some basic properties of MN-series. We give
the residue formula in section 3. Then we talk about the (diagonal
and non-diagonal) Lagrange inversion formulas in section 4, and
give a new proof and a generalization of Dyson's conjecture in
section 5.

\section{Basic Properties of Malcev-Neumann Series}

A \emph{totally ordered abelian group} or TOA-group is an abelian
group $\cal{\mathcal{G}}$ (written additively) equipped with a
total ordering $\le$ that is compatible with the group structure
of $\mathcal{G}$; i.e., for all $x,y,z\in \mathcal{G}$, $x<y$
implies $x+z<y+z$. Such an ordering $<$ is also called {\em
translation invariant}. The abelian groups $\ZZ,$ $\QQ,$ and $\RR$
are all totally ordered abelian groups under the natural ordering.

Let $K$ be a field. A formal series $\eta$ on $\mathcal{G}$ over
$K$ has the form
$$\eta =\sum_{g\in \mathcal{G}} a_g t^g, $$
where $a_g\in K$ and $t^g$ is regarded as a symbol. Let
$\tau=\sum_{h\in \mathcal{G}} b_h t^h$ be another formal series on
$\mathcal{G}$. Then the \emph{product} $\eta \tau$ is defined if
for every $f\in \mathcal{G}$, there are only finitely many pairs
$(g,h)$ of elements of $\mathcal{G}$ such that $a_g$ and $b_h$ are
nonzero and $g+h=f$. In this case,
$$\eta \tau :=\sum_{f\in \mathcal{G}} t^f\sum_{g+h=f}a_gb_h.$$
The \emph{support} $\supp (\eta)$ of $\eta$ is defined to be $\{\,
g\in \mathcal{G}: a_g\ne 0\, \}$.

{}For a TOA-group $\mathcal{G}$, a {\em Malcev-Neumann series}
(MN-series for short) is a formal series on $\mathcal{G}$ that has
a well-ordered support. Recall that a well-ordered set is a
totally ordered set such that every nonempty subset has a minimum.
We define $K_w[\mathcal{G}]$ to be the set of all such MN-series.

By a theorem of Malcev and Neumann
\cite[Theorem~13.2.11]{passmann}, $K_w[\mathcal{G}]$ is a field
for any TOA-group. A sketch of the proof will be introduced since
we will use some of the facts later.

Let us see some examples of MN-series first.
\begin{enumerate}
    \item $K_w[\ZZ]\simeq K((x))$ is the field of Laurent series.
    \item $K_w[\QQ]$ strictly contains the field $K^\fra ((x))$ of fractional Laurent
series \cite[p.~161]{EC2}, and is more complicated. It includes as
a subfield the generalized Puiseux field \cite{vaidya} with
respect to a prime number $p$, which  consists all series $f(x)$
such that $\supp(f)$ is a well-ordered subset of $\QQ$ and there
is an $m$ such that for any $\alpha \in \supp(f)$ we have
$m\alpha=n_\alpha /p^{i_\alpha}$ {for some integer $n_\alpha$ and
nonnegative integer $i_\alpha$}.
\item Let $\QQ^\times$ be the multiplicative group of positive rational numbers. Then $\QQ^\times$ is
a TOA-group, and $K_w[\QQ^\times]$ is a field of MN-series.
\end{enumerate}

The set of MN-series $K_w[\mathcal{G}]$ is clearly closed under
addition. The following proposition is the key to showing that
$K_w[\mathcal{G}]$ is closed under multiplication, so that
$K_w[\mathcal{G}]$ is a ring.

{}For two subsets $A$ and $B$ of $G$, we denote by $A+B$ the set
$\{\, a+b \mid a\in A,b\in B\,\}$.
\begin{prop}[\cite{passmann}, Lemma~13.2.9]
\label{p-welplus} If $\mathcal{G}$ is a TOA-group and $A,B$ are
two well-ordered subsets of $\mathcal{G}$ then $A+B$ is also
well-ordered.
\end{prop}

{}For a TOA-group $\mathcal{G}$, $K_w[\mathcal{G}]$ is a maximal
ring in the set of all formal series on $\mathcal{G}$: if $\eta
=\sum_{g\in \mathcal{G}} a_g t^g$ is not in $K_w[\mathcal{G}]$,
then adding $\eta $ into $K_w[\mathcal{G}]$ cannot yield a ring.
{}For if $\supp (\eta)$ is not well-ordered, we can assume that
$g_1>g_2>\cdots $ is an infinite decreasing sequence in
$\supp(\eta)$. Let $\tau= \sum_{n\ge 1} a_{g_n}^{-1} t^{-g_n}$.
Note that $\tau \in K_w[\mathcal{G}]$, since $-g_1<-g_2<\cdots $
is well-ordered. But the constant term of $\eta \tau$ equals an
infinite sum of $1$'s, which diverges.

Let $[t^g] \eta$ be the coefficient of $t^g$ in $\eta$. Let
$\eta_1,\eta_2,\dots $ be a series of elements in
$K_w[\mathcal{G}]$. Then we say that $\eta_1+\eta_2+\cdots $ {\em
strictly converges} to $\eta \in K_w[\mathcal{G}]$, if for every
$g\in \mathcal{G}$, there are only finitely many $i$ such that
$[t^g] \eta_i\ne 0$, and $\sum_{i\ge 1} [t^g] \eta_i = [t^g]
\eta$. If $\eta_1+\eta_2+\cdots $ strictly converges to some $\eta
\in K_w[\mathcal{G}]$, then we say that $\eta_1+\eta_2+\cdots$
\emph{exists} (in $K_w[\mathcal{G}]$). Note that $\sum_{n\ge 1}
2^{-n}$ does not strictly converge to $1$.

Let $f(z)=\sum_{n\ge 0} b_n z^n$ be a formal power series in
$K[[z]]$, and let $\eta\in K_w[\mathcal{G}]$. Then we define the
composition $f\circ \eta$ to be
$$f\circ \eta := f(\eta)= \sum_{n\ge 0} b_n \eta^n $$
if the sum exists.

If $\eta\ne 0$ belongs to $K_w[\mathcal{G}]$, then it has a
nonempty well-ordered support so that we can define the {\em
order} of $\eta$ to be $\ord(\eta)=\min \left(\supp(\eta)\right)$.
The {\em initial term} of $\eta$ is the term with the smallest
order. It is clear that $\ord(\eta \tau)=\ord(\eta)+\ord(\tau)$.
The order of $0$ is treated as $\infty$.

\begin{thm}[Composition Law]
\label{t-welcomposition} If $f\in K[[z]]$ and $\eta \in
K_w[\mathcal{G}]$ with $\ord (\eta)>0$, then $f\circ \eta$
strictly converges in $K_w[\mathcal{G}]$.
\end{thm}

The detailed proof of this composition law can be found in
\cite[Chapter~3.1]{xinthesis}. It consists of two parts: one is to
show that for any $g\in \mathcal{G}$, $[t^g] f\circ \eta$ is a
finite sum of elements in $K$; the other is to show that the
support of $f\circ \eta$ is well-ordered. The following
proposition is the key to the proof.

We denote by $A^{+n}$ the set $A+A+\cdots +A$ of $n$ copies of
$A$. A subset $A$ of $\mathcal{G}$ is said to be {\em positive},
denoted by $A>0$, if $a>0$ for all $a\in A$.
\begin{prop}[\cite{passmann}, Lemma 13.2.10]\label{p-welpowerunion}
Let $\mathcal{G}$ be a TOA-group.
If $A$ is a positive well-ordered subset of $\mathcal{G}$, then
$\cup_{n\ge 0} A^{+ n}$ is also well-ordered.
\end{prop}

\begin{cor}\label{c-fieldstructure}
{}For any $\eta\in K_w[\mathcal{G}]$ with initial term $1$,
$\eta^{-1}\in K_w[\mathcal{G}]$.
\end{cor}
\begin{proof}
Write $\eta=1-\tau$. Then $\tau\in K_w[\mathcal{G}]$ and $\ord(\tau)>0$. By Theorem 
\ref{t-welcomposition}, $\sum_{n\ge 0} \tau^n$ strictly converges
in $K_w[\mathcal{G}]$. Knowing that $[t^g] (1-\tau) \cdot
\sum_{n\ge 0} \tau^n$ is a finite sum for every $g$, we can check
that $(1-\tau) \cdot \sum_{n\ge 0} \tau^n$ reduces to $ 1$ after
cancelation.
\end{proof}

So for any $\eta\in K_w[\mathcal{G}]$ with initial term $f$,
$\eta=f(1-\tau)$ with $\ord(\tau)>0$, and the expansion of
$\eta^{-1}$ is given by $f^{-1}\sum_{n\ge 0}\tau^n$. This implies
that $K_w[\mathcal{G}]$ is a field.

\begin{dfn}
If $\mathcal{G}$ and $\mathcal{H}$ are two TOA-groups, then the
\emph{Cartesian product} $\mathcal{G}\times \mathcal{H}$ is
defined to be the set $\mathcal{G}\times \mathcal{H}$ equipped
with the usual addition and the reverse lexicographic order, i.e.,
$(x_1,y_1)\le (x_2,y_2)$ if and only if $y_1<_\mathcal{H} y_2$ or
$y_1=y_2$ and $x_1\le_\mathcal{G} x_2$.
\end{dfn}

We define $\mathcal{G}^n$ to be the Cartesian product of $n$
copies of $\mathcal{G}$. It is an easy exercise to show the
following.

\begin{prop}
The Cartesian product of finitely many TOA-groups is a TOA-group.
\end{prop}

One important example is $\ZZ^n$ as a totally ordered abelian
group.

When considering the ring $K_w(\mathcal{G}\times \mathcal{H})$, it
is natural to treat ${(g,h)}$ as $g+h$, where $g$ is identified
with $(g,0)$ and $h$ is identified with $(0,h)$. With this
identification, we have the following.

\begin{prop}\label{p-malcev-prod}
The field $K_w[\mathcal{G}\times \mathcal{H}]$ is the same as the
field $\left(K_w[\mathcal{G}]\right)_w[\mathcal{H}]$ of
Malcev-Neumann series on $\mathcal{H}$ with coefficients in
$K_w[\mathcal{G}]$.
\end{prop}
\begin{proof}
Let $\eta \in K_w[\mathcal{G}\times \mathcal{H}]$, and let
$A=\supp(\eta)$. Let $p$ be the second projection of
$\mathcal{G}\times \mathcal{H}$, i.e., $p(g,h)=h.$

We first show that $p(A)$ is well-ordered. If not, then we have
{an infinite sequence $(g_1,h_1), (g_2,h_2),\dots$ of elements of
$A$ such that  $p(g_1,h_1)>p(g_2,h_2)>\cdots $,} which by
definition becomes $h_1>h_2>\cdots  $. Then in the reverse
lexicographic order, this implies that $(g_1,h_1)>(g_2,h_2)>\cdots
$ is an infinite decreasing sequence in $A$, a contradiction. So
$p(A)$ is well-ordered.

Now $\eta$ can be written as
$$\eta =\sum_{h\in p(A)} \biggl(\sum_{g\in \mathcal{G}, (g,h)\in A}  a_{g,h}
t^g\biggr) t^h.$$ Since for each $h\in p(A)$, the set $\{\, g\in
\mathcal{G}: (g,h)\in A \,\}$ is a clearly a well-ordered subset
of $\mathcal{G}$, $\sum_{g\in \mathcal{G}, (g,h)\in A} a_{g,h}
t^g$ belongs to $K_w[\mathcal{G}]$ for every $h$, and hence
$\eta\in \left(K_w[\mathcal{G}]\right)_w[\mathcal{H}]$.

Conversely, let $\tau=\sum_{h\in D} b_h t^h \in
(K_w[\mathcal{G}])_w[\mathcal{H}]$, where $D=\supp(\tau)$ is a
well-ordered subset of $\mathcal{H}$, and $b_h\in
K_w[\mathcal{G}]$. Let $B_h$ denote the support of $b_h$.
 We need to show that $\bcup_{h\in D} (B_h\times \{\,h \,\})$ is
well-ordered in $\mathcal{G}\times \mathcal{H}$. Let $A$ be any
{nonempty} subset of $\bcup_{h\in D} (B_h\times \{\,h \,\})$. We
show that $A$ has a smallest element. Since $p(A)$ is a subset of
the well-ordered set $D$, we can take $h_0$ to be the smallest
element of $p(A)$. The set $A\cap (B_{h_0}\times \{\,h_0 \,\})$ is
well-ordered for it is a subset of the well-ordered set
$B_{h_0}\times \{\,h_0 \,\}$. Let $(g_0,h_0)$ be the smallest
element of $A\cap (B_{h_0}\times \{\,h_0 \,\})$. Then $(g_0,h_0)$
is also the smallest element of $A$.
\end{proof}

Let $K$ be a field. The field of iterated Laurent series $K\ll
x_1,\dots ,x_n \gg$ is inductively defined to be the field of
Laurent series in $x_n$ with coefficients in $K\ll
x_1,\dots,x_{n-1}\gg$, with $K\ll x_1\gg$ being the field of
Laurent series $K((x_1))$.
\begin{cor}
$$K_w[\ZZ^n] \simeq K\ll x_1,x_2,\dots,x_n\gg.$$
\end{cor}
{The detailed proof of this corollary is left to the reader.} We
only describe the identification as follows. Let $\{\,
e_i\,\}_{1\le i\le n}$ be the standard basis of $\ZZ^n$. Then
$x_i$ is identified with $t^{e_i}$. The field of iterated Laurent
series turns out to be the most useful special kind of MN-series
\cite{xiniterate,xinthesis}.

We conclude this section with the following remark.
\begin{rem}
MN-series were originally defined on totally ordered groups. It
was shown in \cite[Chapter~3.1]{xinthesis} that the results in
this section
{can be generalized:} $\mathcal{G}$ can be replaced with a totally
ordered monoid (a semigroup with a unit), and $K$ can be replaced
with a commutative ring with a unit.
\end{rem}

\section{The Residue Theorem}

Observe that any subgroup of a TOA-group is still a TOA-group
under the induced total ordering. Let $\mathcal{G}$ be a TOA-group
and let $\mathcal{H}$ be {an abelian} group. If
$\rho:\mathcal{H}\to \mathcal{G}$ is an injective homomorphism,
then $\rho(\mathcal{H})\eql \mathcal{H}$ is a subgroup of
$\mathcal{G}$. We can thus regard $\mathcal{H}$ as a subgroup of
$\mathcal{G}$ through $\rho$. The induced ordering $\le^\rho$ on
$\mathcal{H}$ is given by $h_1\le^\rho h_2 \Leftrightarrow
\rho(h_1)\le_\mathcal{G} \rho(h_2)$. Thus $\mathcal{H}$ is a
TOA-group under $\le^\rho$. Clearly a subset $A$ of
$(\mathcal{H},\le^\rho)$ is well-ordered if and only if $\rho(A)$
is well-ordered in $(\mathcal{G},\le_\mathcal{G})$.

Let $\mathcal{G}$ be a TOA-group. We can give $\mathcal{G}$ a
different ordering so that under this new ordering $\mathcal{G}$
is still a TOA-group. {}For instance, the total ordering ${\le}^*$
defined by $g_1\le g_2 \Leftrightarrow g_2\le^* g_1$ is clearly
such an ordering. One special class of total orderings is
interesting for our purpose. If $\rho:\mathcal{G}\to \mathcal{G}$
is an injective endomorphism, then the induced ordering $\le^\rho$
is also a total ordering on $\mathcal{G}$. We denote the
corresponding field of MN-series by $K_w^\rho[\mathcal{G}]$.

{}For example, if $\mathcal{G}=\ZZ^n$, then any nonsingular matrix
$M\in GL(\ZZ^n)$ induces an injective endomorphism. In particular,
$K_w[\ZZ^2]\eql K\ll x,t \gg$ is the field of double Laurent
series, and $K_w^\rho[\ZZ^2]\eql K\ll \xx,t\gg$, where the matrix
corresponding to $\rho$ is the diagonal matrix $\diag (-1,1)$. It
is easy to see that $K\ll x_1^{\epsilon_1},\dots
,x_n^{\epsilon_n}\gg$ with $\epsilon_i=\pm 1$ are special fields
of MN-series $K^\rho\ll x_1,\dots ,x_n \gg$, where the
corresponding matrix for $\rho$ is the diagonal matrix with
entries $\epsilon_i$.

Series expansions in a field of MN-series depend on the total
ordering $\le^\rho$. When comparing monomials, it is convenient to
use  the symbol $\preceq^\rho$: {if $g_1 \le^\rho g_2$ then we
write $t^{g_1}\preceq^\rho t^{g_2}$.} We shall call attention to
the expansions in the following example.

Let $\rho$ be defined by $\rho(x)=x^2y$ and $\rho(y)=xy^2$, and
consider $K^\rho\ll x,y\gg$. The expansion of $1/(x-y)$ is given
by $$\frac{1}{x-y}=\frac{1}{x}\cdot
\frac{1}{1-y/x}=\frac{1}{x}\sum_{k\ge 0} y^k/x^k,$$ since
$\rho(y/x)=\rho(y)/\rho(x)=y/x \succ 1$, which implies $1
\prec^\rho y/x$.

Now notice the expansion of $1/(x^2-y)$ is given by
$$\frac{1}{x^2-y}=-\frac{1}{y}\cdot \frac{1}{1-x^2/y}=-\frac{1}{y}\sum_{k\ge 0} x^{2k}/y^k,$$ since
$\rho (y/x^2) =\rho(y)/\rho(x^2)=1/x^3 \prec 1$, which implies
$1\prec^\rho x^2/y$.

In order to state the residue theorem, we need more concepts.
Consider the following situation. Let $\mathcal{G}$ and
$\mathcal{H}$ be groups with $\mathcal{H}\eql \ZZ^n$, and suppose
that we have a total ordering $\le $ on the direct sum
$\mathcal{G}\oplus \mathcal{H}$ such that $\mathcal{G}\oplus
\mathcal{H}$ is a TOA-group. We identify $\mathcal{G}$ with
$\mathcal{G}\oplus 0$ and $\mathcal{H}$ with $0\oplus
\mathcal{H}$. Let $e_1,e_2,\dots ,e_n$ be a basis of
$\mathcal{H}$. Let $\rho$ be the endomorphism on
$\mathcal{G}\oplus \mathcal{H}$ that is generated by $\rho(e_i)=
g_i+\sum_j m_{ij} e_j$ for all $i$, where $g_i\in \mathcal{G}$,
and $\rho(g)=g$ for all $g\in \mathcal{G}$. Then $\rho$ is
injective if the matrix $M=\left(m_{ij}\right)_{1\le i,j\le n}$
belongs to $GL(\ZZ^n)$, i.e., $\det(M)\ne 0$.

It is natural to use new variables $x_i$ to denote $t^{e_i}$ for
all $i$. Thus monomials in $K_w[\mathcal{G}\oplus \mathcal{H}]$
can be represented as $ t^gx_1^{k_1}\cdots x_n^{k_n}$.
Correspondingly, $\rho$ acts on monomials by $\rho( t^g)=t^g$ for
all $g\in \mathcal{G}$, and $\rho(x_i)= t^{g_i} x_1^{m_{i1}}\cdots
x_n^{m_{in}}.$

\noindent {\bf Notation}: If $f_i$ are monomials, we use $\mb{f}$
to denote the homomorphism $\rho$ generated by $\rho(x_i)=f_i$.

 \vspace{3mm} An
element $\eta\ne 0$ of $K_w[\mathcal{G}\oplus \mathcal{H}]$ can be
written as
$$\eta=\sum_{\mb{k}\in \ZZ^n} \sum_{g\in \mathcal{G}} a_{g,\mb{k}} t^g x_1^{k_1}\cdots x_n^{k_n}
= \sum_{\mb{k}\in \ZZ^n} b_{\mb{k}} \mb{x}^{\mb{k}},$$ where
$a_{g,\mb{k}}\in K$ and $b_{\mb{k}}\in K_w[\mathcal{G}]$. If
$b_{\mb{k}}\mb{x}^{\mb{k}}\ne 0$, then we call it an
\emph{$x$-term} of $\eta$. Since the set $\{\,
\ord(b_{\mb{k}}\mb{x}^{\mb{k}}): \mb{k}\in \ZZ^n, b_{\mb{k}}\ne
0\}$ is a nonempty subset of $\supp(\eta)$, it is well-ordered and
hence has a least element. Because of the different exponents in
the $x$'s, no two of $\ord(b_{\mb{k}}\mb{x^k})$ are equal. So we
can define the $x$-{\em initial} term of $\eta$ to be the $x$-term
that has the least order.

To define the operators $\frac{\partial}{\partial x_i}$,
$\ct_{x_i}$, $\res_{x_i}$, it suffices to consider the case
$\mathcal{H}=\ZZ$. These operators are defined by:
\begin{align*}
\frac{\partial}{\partial x} \sum_{n\in \ZZ} b_n x^n = \sum_{n\in
\ZZ} n b_n x^{n-1}, \qquad \qquad \ct_x \sum_{n\in \ZZ} b_n x^n =
b_0,\qquad\qquad \res_x \sum_{n\in \ZZ} b_n x^n = b_{-1}.
\end{align*}
Multivariate operators are defined by iteration. All these
operators work nicely in the field of MN-series
$K_w[\mathcal{G}\oplus \mathcal{H}]$, because an MN-series has a
well-ordered support, and still has a well-ordered support after
applying these operators.

There are several computational rules
\cite[Lemma~3.2.1]{xinthesis} for evaluating constant terms in the
univariate case, but we are going to concentrate on the residue
theorem in the multivariate case.

\vspace{3mm} In what follows, we suppose $F_i\in
K_w[\mathcal{G}\oplus \mathcal{H}]$ for all $i$.

 \begin{dfn}
{The Jacobian determinant (or simply Jacobian) of
$\mb{F}=(F_1,F_2,\dots, F_n)$ with respect to $\mb{x}$ is defined
to be}
$$J\left({\mb{F}}|{\mb{x}}\right):=
J\left(\frac{F_1,F_2,\ldots ,F_n}{x_1,x_2,\dots ,x_n}\right)
=\det\left(\pad{x_j}{F_i} \right)_{1\le i,j\le n}.$$
\end{dfn}
 When the $x$'s  are clear, we write $J(F_1,F_2,\dots
,F_n)$ for short.

\begin{dfn}
If the $x$-initial term of $F_i$ is $a_i x_1^{b_{i1}}\cdots
x_n^{b_{in}}$, then the Jacobian number of $\mb{F}$ with respect
to $\mathbf{x}$ is defined to be
$$j\left({\mb{F}}|{\mb{x}}\right):=
j\left(\frac{F_1,F_2,\ldots ,F_n}{x_1,x_2,\dots ,x_n}\right)
=\det\left(b_{ij} \right)_{1\le i,j\le n}.$$
\end{dfn}

\begin{dfn}
The log Jacobian of $F_1,\dots ,F_n$ is defined to be
$$LJ(F_1,\dots,F_n):=\frac{x_1\cdots x_n}{F_1 \cdots F_n} J(F_1,\dots ,F_n).$$
\end{dfn}
We call it the log Jacobian because formally it can be written as
\citep{wilson}
$$LJ(F_1,\dots,F_n)=J\left(\frac{\log F_1,\ldots ,\log F_n}{\log x_1,
\dots ,\log x_n}\right),$$ since
$$\frac{\partial \log F}{\partial \log x}= \frac{\partial \log F}{\partial F}
\frac{\partial F}{\partial \log x}=\frac{1}{F}\frac{\partial
{}F}{\partial  x} \frac{\partial x}{\partial \log
x}=\frac{x}{F}\frac{\partial F}{\partial x}.$$

\begin{rem}
The Jacobian is convenient in residue evaluation, while the log
Jacobian is convenient in constant term evaluation.
\end{rem}

The following lemma is needed for the proof of our residue
theorem. It is also a kind of generalized composition law.

Let $\Phi$ be a formal series in $x_1,\dots, x_n$ with
coefficients in $K_w[\mathcal{G}]$, and let $F_i\in
K_w[\mathcal{G}\oplus \mathcal{H}]$. Then $\Phi(F_1,\dots ,F_n)$
is obtained from $\Phi$ by replacing $x_i$ with $F_i$. The
following lemma gives a simple sufficient condition for the
convergence of $\Phi(F_1,\dots, F_n)$.

\begin{lem}
\label{l-malcev-comp} Let $\Phi$ and $F_i$ be as above and let
$f_i$ be the initial term of $F_i$ for all $i$. Suppose
$j(F_1,\dots ,F_n)\ne 0$. Then $\Phi(x_1,\dots ,x_n)\in K_w^\mb{f}
[\mathcal{G}\oplus \mathcal{H}]$ if and only if
$\Phi(f_1,\dots,f_n)$ exists in $ K_w[\mathcal{G}\oplus
\mathcal{H}]$, and if these conditions hold then $\Phi(F_1,\dots
{}F_n)$ exists in $ K_w[\mathcal{G}\oplus \mathcal{H}]$.
\end{lem}

\begin{proof}[Proof of Lemma \ref{l-malcev-comp}]
 We first show the equivalence.
The map $\rho: x_i\to {f_i}$ induces an endomorphism on
$\mathcal{H}\simeq \ZZ^n$. This endomorphism is injective since
$j(f_1,\dots ,f_n)\ne 0$, which is equivalent to $j(F_1,\dots
,{}F_n)\ne 0$. Therefore $\rho$ also induces an injective
endomorphism on $\mathcal{G}\oplus \mathcal{H}$. We see that
 $\supp(\Phi(f_1,\dots, f_n))$  is well-ordered in $\mathcal{G}\oplus \mathcal{H}$ if and only
if $\rho \left(\supp(\Phi(x_1,\dots, x_n))\right)$ is
well-ordered. This, by definition, is to say that $
\Phi(x_1,\dots, x_n)\in K_w^\mb{f} [\mathcal{G}\oplus
\mathcal{H}].$

Now we show the implication. Write each $F_i$ as $f_i(1+\tau_i)$,
with $\ord(\tau_i)>0$. Given the convergence of
$\Phi(f_1,\dots,f_n)$ we first show that for every $g\in
\mathcal{G}$ and $\mb{m}\in \ZZ$, $[t^g \mb{x}^{\mb{m}}] \,
\Phi(F_1,\dots ,F_n)$ is a finite sum.

Write $\Phi$ as $\sum_{\mb{k}\in \ZZ^n} a_{\mb{k}} \mb{x}^\mb{k}$.
Let $A$ be the support of $\Phi(\mb{f})$. Then $A$ is the disjoint
union of $\supp (a_{\mb{k}} \mb{f}^\mb{k})$ for all $\mb{k}$. This
follows from the first part: $\rho$ is injective.

Now \begin{equation} \label{e-3composition-law} \Phi(F_1,\dots
,F_n)=\sum_{\mb{k}\in \ZZ^n} a_{\mb{k}}
\mb{f}^\mb{k}(1+\tau_1)^{k_1}\cdots (1+\tau_n)^{k_n}.
\end{equation}

We observe that replacing any nonzero element in $K$ by $1$ will
not reduce the number of summands, so $(1+\tau_i)^{k_i}$ can be
replaced with $(1-\tau_i)^{-1}=\sum_{l\ge 0} \tau_i^{l}$.
Therefore, the number of summands for the coefficient of $t^g
\mb{x}^{\mb{m}}$ in $\Phi(F_1,\dots,F_n)$ is no more than that in
$$ \sum_{\mb{k}\in \ZZ^n} a_{\mb{k}}
\mb{f}^\mb{k}(1-\tau_1)^{-1}\cdots
(1-\tau_n)^{-1}=(1-\tau_1)^{-1}\cdots
(1-\tau_n)^{-1}\sum_{\mb{k}\in \ZZ^n} a_{\mb{k}} \mb{f}^\mb{k},$$
which is a finite product of elements in $K_w[\mathcal{G}\oplus
\mathcal{H}]$. Note that in obtaining the right-hand side of the
above equation, we used the fact that the supports of $a_{\mb{k}}
\mb{f}^\mb{k} $ are disjoint for all $\mb{k}$.

The proof of the lemma will be finished after we show that
$\Phi(F_1,\dots,F_n)$ has a well-ordered support. Let $T_i$ be the
support of $\tau_i$. Then the support of $(1+\tau_i)^{k_i}$ is
contained in $\cup_{l\ge 0} T_i^{+l}$. Thus for every $\mb{k}$
$$\supp a_{\mb{k}}
\mb{f}^\mb{k}(1+\tau_1)^{k_1}\cdots  (1+\tau_n)^{k_n} \subseteq
A+\cup_{l\ge 0} T_1^{+l}+\dots +\cup_{l\ge 0} T_n^{+l},
$$
which is well-ordered by Proposition \ref{p-welplus} and
Proposition \ref{p-welpowerunion}. So by
\eqref{e-3composition-law}, the support of $\Phi(F_1,\dots,F_n)$
is also well-ordered.

\end{proof}

\begin{rem}
The implication in Lemma \ref{l-malcev-comp} is not true when
$j(F_1,\dots,F_n)=0$. For instance, let $\Phi=\sum_{k\ge 0}
x_2^k/x_1^k -\sum_{k\ge 0} x_2^{3k}/x_1^{2k}$ and let $F_1=x_1^2$,
$F_2=x_1(1+x_1)$. Then it is straightforward to check that
$\Phi(f_1,f_2)=0$, but $\Phi(F_1,F_2)$ is not in $K\ll x_1\gg$.
\end{rem}

\vspace{3mm} \noindent {\bf Notation.} Starting with a TOA-group
$\mathcal{G}\oplus \mathcal{H}$ as described above, let $\Phi$ be
a formal series on $\mathcal{G}\oplus \mathcal{H}$. When we write
$\ct_\mathbf{x}^\rho \Phi(x_1,\dots ,x_n)$, we mean both that
$\Phi(x_1,\dots ,x_n)$ belongs to $ K_w^\rho[\mathcal{G}\oplus
\mathcal{H}]$, and that the constant term is taken in this field.
When $\rho$ is the identity map, it is omitted. When we write
$\ct_\mb{F} \Phi(F_1,\dots ,F_n)$, it is assumed that
$\Phi(x_1,\dots ,x_n)\in K_w^\mb{f}[\mathcal{G}\oplus
\mathcal{H}]$, {where $f_i$ is the initial term of $F_i$,} and we
are taking the constant term of $\Phi(x_1,\dots ,x_n) $ in the
ring $K_w^\mb{f}[\mathcal{G}\oplus \mathcal{H}]$. Or equivalently,
we always have
$$\ct_\mb{F} \Phi(F_1,\dots ,F_n) = \ct_{\mb{x}}\mbox{}^\mb{f} \Phi(x_1,\dots ,x_n).$$
This treatment is particularly useful  when dealing with rational
functions.

\vspace{3mm} Now comes our residue theorem for
$K_w[\mathcal{G}\oplus \mathcal{H}]$, in which we will see how an
element in one field is related to an element in another field
through taking constant terms.
\begin{thm}[Residue Theorem]\label{t-MNresidue}
Suppose for each $i$, $F_i\in K_w[\mathcal{G}\oplus \mathcal{H}]$
has $x$-initial term $f_i=a_i x_1^{b_{i1}}\cdots x_n^{b_{in}}$
with $a_i\in K_w[\mathcal{G}]$. If $j(F_1,\dots ,F_n)\ne 0$, then
 for any
 $\Phi(\mb{x})\in K_w^\mb{f} [\mathcal{G}\oplus \mathcal{H}]$, we have
\begin{align}\label{e-wel-residue}
\res_{\mb{x}} \Phi(F_1,\dots, F_n) J(F_1,\dots
,{}F_n)=j(F_1,\dots,F_n) \res_{\mb{F}} \Phi(F_1,\dots,F_n).
\end{align}
Equivalently, in terms of constant terms, we have
$$
\ct_\mathbf{x} \Phi(F_1,\dots ,F_n) LJ(F_1,\dots ,F_n)
=j(F_1,\dots ,F_n)\ct_\mb{F} \Phi(F_1,\dots ,F_n)
.\eqno{(\ref{e-wel-residue}')}
$$
\end{thm}
\begin{proof}[Proof of Theorem \ref{t-MNresidue}]
Replace $\Phi(F_1,\dots, F_n)$ with $F_1\cdots F_n \Phi(F_1,\dots
,F_n)$ in \eqref{e-wel-residue}. Then by a
 straightforward algebraic manipulation, we will get (\ref{e-wel-residue}$'$).
Similarly we can obtain (\ref{e-wel-residue}) from
(\ref{e-wel-residue}$'$). This shows the equivalence.

By the hypothesis and Lemma \ref{l-malcev-comp}, the left-hand
side of \eqref{e-wel-residue} exists by taking the constant term
in $K_w[\mathcal{G}\oplus \mathcal{H}]$, while the right-hand side
exists by taking the constant term in $K_w^\mb{f}
[\mathcal{G}\oplus \mathcal{H}]$.

{}For the remaining part it suffices to show that the theorem is
true for monomials $\Phi$ by multilinearity. The proof will be
completed after we show Lemmas \ref{l-3-residue-n} and
\ref{l-3-residue-y} below.
\end{proof}

\begin{rem}
When $j(F_1,\dots ,F_n)= 0$, $\Phi(F_1,\dots ,F_n)$ is only well
defined in some special cases. In such cases, \eqref{e-wel-residue}
also holds. {}For example,
if $\Phi(x_1,\dots ,x_n)$ is a Laurent polynomial, then
$\Phi(F_1,\dots ,F_n)$ always exists.
\end{rem}

\begin{rem}
The theorem holds for any rational function $\Phi$, i.e.,
$\Phi(x_1,\dots,x_n)$ belongs to the quotient field of
$(K_w[\mathcal{G}])[\mathcal{H}]$. This follows from the fact that
$K_w^{\mb{f}}[\mathcal{G}\oplus \mathcal{H}]$ is a field
containing $(K_w[\mathcal{G}])[\mathcal{H}]$ as a subring.
\end{rem}

The proof of our residue theorem and lemmas basically comes from
\citep{reversion}, except for the proof of Lemma
\ref{l-3-residue-y}, which uses the original idea of Jacobi.

The following properties of Jacobians can be easily checked.
\begin{lem}\label{l-lres}
We have
\begin{enumerate}
\item $J(F_1,F_2,\ldots ,F_n)$ is $K_w[\mathcal{G}]$-multilinear.
\item $J(F_1,F_2,\ldots ,F_n)$ is alternating; i.e.,
$J(F_1,F_2,\ldots ,F_n)=0$ if $F_i=F_j$ for some $i\ne j$.
\item $J(F_1,F_2,\ldots ,F_n)$ is anticommutative; i.e.,
$$J(F_1,\ldots,F_i,\ldots ,F_j,\ldots,F_n)=-J(F_1,\ldots,F_j,\ldots ,F_i,\ldots,F_n).$$
\item $($Composition Rule$)$ If $g(z)\in K((z))$ is a series in one variable, then
$$J(g( F_1),F_2,\ldots ,F_n)=\frac{dg}{dz}(F_1) J(F_1,F_2,\ldots ,F_n).$$
\item $($Product Rule$)$ $$J(F_1G_1,F_2,\ldots ,F_n)=F_1J(G_1,F_2,\ldots ,F_n)+
G_1J(F_1,F_2,\ldots ,F_n).$$
\item $J(F_2^{-1},F_2,\ldots ,F_n)=0$.
\end{enumerate}
\end{lem}

A formal series on $\mathcal{G}\oplus \mathcal{H}$ having only one
$x$-term is called an $x$-\emph{monomial}.
\begin{lem} \label{l-residue-f} If all $f_i$ are $x$-monomials, then
\begin{align}\label{e-3-jacobian}
LJ(f_1\dots ,f_n)=j(f_1,\dots ,f_n).\qquad\qquad\
\end{align}
Equivalently,
$$J(f_1,\dots ,f_n)= j(f_1,\dots ,f_n)\frac{f_1\cdots f_n}{x_1\cdots x_n}.\eqno{(\ref{e-3-jacobian}')}$$
\end{lem}
\begin{proof}
Suppose that  for every $i$, $f_i=a_i x_1^{b_{i1}}\cdots
x_n^{b_{in}}$, where $a_i$ is in $K_w[\mathcal{G}]$. Then
$\partial f_i/\partial x_j = b_{ij} f_i/x_j$. {}Factoring $f_i$ from
the $i$th row of the Jacobian matrix for all $i$ and then
factoring $x_j^{-1}$ from the $j$th column for all $j$, we get
$$J(f_1,f_2,\dots ,f_n)=\frac{f_1\cdots f_n}{x_1\cdots x_n}\det(b_{ij}).$$
Equation \eqref{e-3-jacobian} and (\ref{e-3-jacobian}$'$) are just
rewriting of the above equation.
\end{proof}

\begin{lem}\label{l-residue}
$$\res_\mathbf{x} J(F_1,\ldots , F_n)=0.$$
\end{lem}
\begin{proof}
By multilinearity, it suffices to check $x$-monomials $F_i$.
Suppose $F_i=f_i$ as given in the proof of Lemma
\ref{l-residue-f}. Then equation (\ref{e-3-jacobian}$'$) can be
rewritten as
$$J(F_1,\dots ,F_n)= \det(b_{ij}) a_1\cdots a_n  x_1^{-1+\sum b_{i1}} \cdots x_n^{-1+\sum
{b_{in}}}. $$ If $\sum {b_{i1}}=\sum {b_{i2}}=\cdots =\sum
{b_{in}}=0$, then the Jacobian number is $0$, and therefore the
residue is $0$. Otherwise, at least one of the $x_i$'s has
exponent $\ne -1$,
 so the residue is $0$ by definition.
\end{proof}

\begin{lem}\label{l-3-residue-n}
{}For all integers $e_i$ with at least one $e_i\ne -1$, we have
\begin{align}\label{e-residueth}
\res_\mathbf{x} F_1^{e_1}\cdots F_n^{e_n} J(F_1,\ldots ,F_n)= 0.
\end{align}
\end{lem}
\begin{proof}
The clever proof in \citep[Theorem 1.4]{reversion} also works
here.

Permuting the $F_i$ and using $(3)$ of Lemma \ref{l-lres}, we may
assume that $e_1\ne -1$,\dots, $e_j\ne -1$, but $e_{j+1}=\cdots
=e_{n}=-1$, for some $j$ with $1\le j\le n$. Setting
$G_i=\frac{1}{e_i+1}F_i^{e_i+1}$ for $i=1,\ldots , j$, we have
$$F_1^{e_1}F_2^{e_2}\cdots F_n^{e_n} J(F_1,F_2,\ldots ,F_n)=
{}F_{j+1}^{-1}\cdots F_n^{-1} J(G_1,\ldots, G_j,F_{j+1}, \ldots
,F_n).$$ Then applying the formula
$$F_{j+1}^{-1} J(G_1,\ldots, G_j,F_{j+1}, \ldots ,F_n)=
J(F_{j+1}^{-1}G_1,G_2\ldots, G_j,F_{j+1}, \ldots ,F_n)$$
repeatedly for $j+1,j+2,\dots ,n$, we get
$$ J(F_{j+1}^{-1}\cdots F_n^{-1}G_1,G_2, \ldots, G_j,F_{j+1}, \ldots ,F_n).$$
The result now follows from Lemma \ref{l-residue}.
\end{proof}

{}For the case $e_1=e_2=\cdots =e_n=-1$, we have
\begin{lem}\label{l-3-residue-y}
\begin{align}
\label{e-3-residue-y}
 \res_{\mb{x}}F_1^{-1}\cdots F_n^{-1}J(F_1,\dots ,F_n)=j(F_1,\dots, F_n).
\end{align}
\end{lem}
The simple proof for this case in \citep{reversion} does not apply
in our situation. The reason will be explained in Proposition
\ref{p-3-ljacobian}.

Note that Lemma \ref{l-3-residue-y} is equivalent to saying that
\begin{align}
\label{e-3-lresidue-y}
 \ct_{\mb{x}}LJ(F_1,\dots ,F_n)=j(F_1,\dots, F_n).
\end{align}

\begin{proof}
Let $f_i:=a_ix_1^{b_{i1}}\cdots x_n^{b_{in}}$ be the $x$-initial
term of $F_i$. Then $F_i=f_i B_i$, where $B_i\in
K_w[\mathcal{G}\oplus \mathcal{H}]$ has $x$-initial term $1$. By
the composition law,
 $\log ( B_i) \in
K_w[\mathcal{G}\oplus \mathcal{H}] $. Now applying the product
rule, we have
\begin{align*}
&F_1^{-1}\cdots F_n^{-1}J(F_1,F_2,\ldots ,F_n) \\
&\qquad \qquad= f_1^{-1}F_2^{-1}\cdots F_n^{-1} J(f_1,F_2,\ldots,
{}F_n)
+B_1^{-1}F_2^{-1}\cdots F_n^{-1} J(B_1,F_2,\ldots, F_n)\\
&\qquad\qquad=f_1^{-1}F_2^{-1}\cdots F_n^{-1} J(f_1,F_2,\ldots,
{}F_n)+ F_2^{-1}\cdots F_n^{-1} J( \log (B_1),F_2,\ldots, F_n).
\end{align*}
{}From Lemma \ref{l-3-residue-n}, the last term in the above
equation has no contribution to the residue in $x$, and hence can
be discarded.

The same procedure can be applied to $F_2,F_3,\ldots ,F_n$.
{}Finally we will get
$$\res_x F_1^{-1}\cdots F_n^{-1}J(F_1,F_2,\ldots ,F_n)=
\res_x f_1^{-1}\cdots f_n^{-1} J(f_1,f_2,\ldots ,f_n),$$ which is
equal to the Jacobian number by Lemma \ref{l-residue-f}.
\end{proof}

The proof of our residue theorem is now completed.

\medskip The next result gives a good reason for using the log
Jacobian.
\begin{prop}\label{p-3-ljacobian}
The $x$-initial term of the log Jacobian $LJ(F_1,\dots ,F_n)$
equals the Jacobian number $j(F_1,\dots ,F_n)$ when it is nonzero.
\end{prop}
\begin{proof}
{}From the definition,
$$LJ(F_1,\dots ,F_{n})= \frac{x_1\cdots x_n}{F_1\cdots F_n}
J(F_1,\dots, F_n)=\frac{x_1\cdots x_n}{F_1\cdots
{}F_n}\sum_{\mb{g}}J(g_1,\dots,g_n) ,$$ where the sum ranges over
all $x$-terms $g_i$ of $F_i$. Applying Lemma \ref{l-residue-f}
gives us
$$LJ(F_1,\dots,F_n)=\sum_{\mathbf{g}}\frac{g_1\cdots g_n}{F_1\cdots F_n} j(g_1,\dots,g_n).$$
The Jacobian number is always an integer. The displayed summand
has the smallest order when $g_i$ equals the $x$-initial term of
$F_i$ for all $i$. It is clear now that we can write
$$LJ(F_1,\ldots ,F_n )= j(F_1,\dots ,F_n)+ \text{higher ordered terms}.$$
To show that $j(F_1,\dots ,F_n)$ is the $x$-initial term, we need
to show that all the other terms that are independent of $x$
cancel. (Note that we do not have this trouble when all the
coefficients belong to $K$.) This is equivalent to saying that
$$ \ct_{\mb{x}}LJ(F_1,\ldots ,F_n )= j(F_1,\dots ,F_n),$$
which follows from Lemma \ref{l-3-residue-y}.
\end{proof}

\begin{exa}
Let $K\ll x,t \gg$ be the working field. Let $F=x^2+xt+x^3t$. Then
the $x$-initial term of $F$ is $x^2$, so $j(F|x)=2$. Now let us
see what happens to the log Jacobian $LJ(F|x)$ of $F$ with respect
to $x$.
\end{exa}
\begin{align*}
LJ(F|x)=\frac{x}{F}\frac{\partial F}{\partial x}
&=\frac{x(2x+t+3x^2t)}{x^2(1+t/x+xt)} \\
&=(2+t/x+3xt) \sum_{k\ge 0} (-1)^k (t/x+xt)^k
\end{align*}
Since every other monomial is divisible by $t$, the initial term
of $LJ(F|x)$ is $2$. It then follows that the $x$-initial term of
$LJ(F|x)$ must contain $2$ and therefore must be the constant term
in $x$.

It is not clear that $2$ is the unique term in the expansion of
$\ct_x LJ(F|x)$, but all the other terms cancel. We check as
follows.
\begin{align*}
\ct_x LJ(F|x) &= \ct_x (2+t/x+3xt) \sum_{k\ge 0} (-1)^k (t/x+xt)^k \\
&= 2 \sum_{k\ge 0 }\binom{2k}{k}t^{2k}- t\sum_{k\ge
0}\binom{2k+1}{k}t^{2k+1}-3t \sum_{k\ge 0} \binom{2k+1}{k+1} t^{2k+1} \\
&= 2+ \sum_{k\ge 1} \left( 2 \binom{2k}{k}-4
\binom{2k-1}{k}\right)t^{2k}.
\end{align*}
Now it is easy to see that the terms, other than $2$, not
containing $x$ in the expansion of the log Jacobian really cancel.
\qed

\medskip {}From Theorem \ref{t-MNresidue} and Lemma
\ref{l-residue-f}, we see directly the following result.
\begin{cor}\label{c-3-residue-m}
If $f_i$ are all $x$-monomials in $K_w[\mathcal{G}\oplus
\mathcal{H}]$,
{$j(f_1,\dots ,f_n)\ne0$, and  $\Phi\in
K_w^\mb{f}[\mathcal{G}\oplus \mathcal{H}]$, then}
$$\ct_{\mb{x}} \Phi(f_1,\dots ,f_n)=\ct_{f_1,\dots ,f_n} \Phi(f_1,\dots ,f_n).$$
\end{cor}

In the case that all $f_i$ are monomials in $K[\mb{x},\mb{\xx}]$
with $j(\mb{f})\ne 0$, $\Phi$ is in $K[\mb{x},\mb{\xx}]$ if and
only $\Phi(f_1,\dots ,f_n)$ is (with possible fractional
exponents). Since $\Phi$ has a finite support, its series
expansion is independent of the working field. In particular, we
have
$$\ct_{f_1,\dots ,f_n} \Phi(f_1,\dots ,f_n)=\ct_{x_1,\dots ,x_n} \Phi(x_1,\dots ,x_n).$$

More generally, we have the following as a consequence of
Corollary \ref{c-3-residue-m} and the above argument.
\begin{cor}\label{c-3-monomial}
Suppose $\mb{y}$ is another set of variables. If $\Phi \in
K[\mb{x},\mb{\xx}]\ll \mb{y}\gg$, and if $f_i$ are all monomials
in $\mb{x}$ with $j(\mb{f})\ne 0$, then
$$\ct_{\mb{x}} \Phi(f_1,\dots ,f_n)=\ct_{\mb{x}} \Phi(x_1,\dots ,x_n).$$
\end{cor}

The following two examples are illustrative in explaining our
residue theorem.
\begin{exa}
The following identity follows trivially by replacing $x$ with
$x^{-1}$.
\begin{align} \label{e-exa-tr1}
\ct_x \sum_{k\ge 0}
x^{-k} = \ct_x  \sum_{k\ge 0} x^k.
\end{align}
\end{exa}
This identity is not as simple as it might appear at first sight.
It equates the constant terms of two elements belonging to two
different fields; namely, the left-hand side of \eqref{e-exa-tr1}
takes the constant term in $K((x^{-1}))$, while the right-hand
side takes the constant term in $K((x))$.

The above cannot be explained by Jacobi's formula, especially when
we write it in terms of rational functions:
\begin{align}\label{e-exa-tr2}
\ct_{x} \frac{1}{1-x^{-1}} = \ct_x \frac{1}{1-x}.
\end{align}
Now let us explain this identity in two ways: one  using our
residue theorem, and the other  using complex analysis.

Let $f=x^{-1}$. Then the log Jacobian $LJ(f|x)=x/f \cdot \partial
f/\partial x = -1$, and the Jacobian number is also $-1$. Thus
$$\ct_x \frac{1}{1-x} =\ct_x \frac{1}{1-f^{-1}} \cdot - LJ(f|x) =\ct_f \frac{1}{1-f^{-1}}.$$
So the $x$ on the left-hand side of \eqref{e-exa-tr2} is indeed
playing the same role with the variable $f$ defined by $f=x^{-1}$.
Now $f^{-1}\succ 1$ since it is the same as $x\succ 1$, and we
have the correct series expansion.

Now we sketch the idea in complex analysis, and describe the
meaning of Jacobian number in the one variable case. We have
$$\ct_x \frac{1}{1-x} = \frac{1}{2\pi i} \oint_\gamma \frac{1}{z(1-z)} dz, $$
where $\gamma$ is the counter-clockwise circle $|z|=\epsilon$ for
sufficiently small positive
$\epsilon$. We can think {of $\epsilon$ as equal to $x$.}

Now if we make a change of variable by $z= 1/u$,
then after simplifying, we get
$$\frac{1}{2\pi i} \oint_{\gamma'} \frac{-1}{u(1-u^{-1})} du
= \ct_f \frac{1}{1-f^{-1}}, $$ where $\gamma'$, the image of
$\gamma$ under the map $z\mapsto 1/u$, is the clockwise circle
$|u|=1/\epsilon$. The Jacobian number $-1$ comes from the
different orientation of the circle. Similarly, if we are making
the change of variable by $z=u^2$, the new circle will be a double
circle, which is consistent with the fact that the Jacobian number
is $2$. \qed

\begin{exa}
Evaluate the following constant term in $K((x))$.
$$\ct_x \frac{(1-x^{-1})^4}{(x-1)(\pi(1-x^{-1})+(1-x^{-1})^2)}.$$
\end{exa}
\begin{proof}[Solution]
Let $F=1-x^{-1}$. Then $LJ(F|x)=x/F\cdot dF/dx=1/(x-1)$. The
$x$-initial term of $F$ is $x^{-1}$ so that the Jacobian number is
$-1$. Hence by our residue theorem, we have
\begin{align*}
\ct_x \frac{(1-x^{-1})^4}{(x-1)(\pi(1-x^{-1})+(1-x^{-1})^2)}&=
\ct_x
\frac{F^4}{\pi F+F^2} LJ(F|x) \\
&=\ct_F\,(-1)\cdot \frac{F^4}{\pi F+F^2}.
\end{align*}
Now the initial term of $F$ is $x^{-1}$ and the initial term of
$F^2$ is $x^{-2}$ so that $F\succ F^2$. Thus the final solution is
$$\ct_F  \frac{-F^2}{1+\pi F^{-1}} =- \pi^2. $$
\end{proof}
\begin{rem}
Suppose the working field is $K((x))$. If the new variable $F$ has
a positive Jacobian number $j(F|x)$, the {second field as
described} in our residue theorem is also $K((x))$. In this case,
Jacobi's formula also applies. If $j(F|x)$ is a negative number,
then we can choose $F^{-1}$ as the new variable to apply Jacobi's
formula. This is why the two fields phenomenon as in the above two
examples {was not noticed} before.
\end{rem}

%
%

The next example is hard to evaluate without using our residue
theorem.
\begin{exa}
Evaluate the following constant term in $\CC\ll x,y,t\gg$.
\begin{align}\label{exa-3-residue}
\ct_{x,y} {x}^{3}{e^{{ {t}/{xy}}}} \left( 2t-3xy \right) \left(
{x}^{3 }y{e^{{ {t}/{xy}}}}-tx-ty \right) ^{-1} \left( x-y \right)
^{-1}
 \left( -1+{x}^{3}{e^{{ {t}/{xy}}}} \right) ^{-1}.
\end{align}
\end{exa}

\begin{proof}[Solution]
The $x$-variables are $x$ and $y$. Let $F=x^2ye^{{t}/{xy}}$,
$G=xy^2 e^{{t}/{xy}}$. It is straightforward to compute the log
Jacobian and the Jacobian number. We have
$$LJ(F,G|x,y)=3-\frac{2t}{xy}, \text{ and }j(F,G|x,y)=3. $$
We can check that \eqref{exa-3-residue} can be written as
$$\ct_{x,y} \frac{F^3G}{(F^2-(F+G)t)(F-G)(G-F^2)}LJ(F,G|x,y).$$
Thus by the residue theorem, the above constant term equals
\begin{align}\label{exa-3-residue1}
\ct_{F,G} \frac{3F^3G}{(F^2-(F+G)t)(F-G)(G-F^2)}=\ct_{F,G}
\frac{3}{ (1-\frac{(F+G)t}{F^2})(1-\frac{G}{F})(1-\frac{F^2}{G})},
\end{align}
where on the right-hand side of \eqref{exa-3-residue1}, we can
check that $1$ is the initial term of each factor in the
denominator.

At this stage, we can use series expansion to obtain the constant
term. We use the following lemma instead.

\begin{lem}\label{l-3-extra}
Suppose that $\Phi$ contains only nonnegative powers in $x$. Then
$$\ct_x \Phi(x) \cdot \frac{1}{1-u/x} = \Phi(u), $$
where $u$ is independent of $x$ and $u\succ x$.
\end{lem}
This lemma is reduced by linearity to the case when $\Phi(x)=x^k$
for some nonnegative integer $k$, which is trivial.

We take the constant term in $G$ first by applying Lemma
\ref{l-3-extra}.
\begin{align*} \ct_{F,G}
\frac{3}{
(1-\frac{(F+G)t}{F^2})(1-\frac{G}{F})(1-\frac{F^2}{G})}&=\ct_F
\frac{3F^3}{(F^2-(F+F^2)t)(F-F^2)}\\
&=\ct_F \frac{3}{(1-t)(1-F)}\cdot \frac{1}{(1-\frac{t}{(1-t)F})} \\
&=\frac{3}{(1-t)(1-\frac{t}{1-t})},
\end{align*}
where in the last step, we applied Lemma \ref{l-3-extra} again.

After simplification, we finally get
$$
\ct_{x,y} {x}^{3}{e^{{ {t}/{xy}}}} \left( 2t-3xy \right) \left(
{x}^{3 }y{e^{{ {t}/{xy}}}}-tx-ty \right) ^{-1} \left( x-y \right)
^{-1}
 \left( -1+{x}^{3}{e^{{ {t}/{xy}}}} \right) ^{-1}=\frac{3}{1-2t}.$$
 \end{proof}

\section{Another View of Lagrange's Inversion Formula\label{ss-lag}}
Let $F_1,\ldots ,F_n$ be power series in variables $x_1,\ldots
,x_n$ of the form $F_i=x_i+$ higher degree terms,  with
indeterminate coefficients for each $i$. It is known, e.g.,
\cite[Proposition~5, p.~219]{reference}, that
$\mathbf{F}=(F_1,\ldots ,F_n)$ has a unique compositional inverse,
i.e., there exists $\mathbf{G}=(G_1,\ldots , G_n)$ where each
$G_i$ is a power series in $x_1,\ldots ,x_n$ such that
$F_i(G_1,\ldots ,G_n)=x_i$
 and $G_i(F_1,\ldots ,F_n)=x_i$ for all $i$.

Lagrange inversion gives a formula for the $G$'s in terms of the
$F$'s. Such a formula is very useful in combinatorics. A good
summary of this subject can be found in \cite{gessel-res}.

The diagonal (or Good's) Lagrange inversion formula deals with the
diagonal case, in which $x_i$ divides $F_i$ for every $i$, or
equivalently, $F_i=x_i H_i$, where $H_i\in K[[x_1,\dots ,x_n]]$
with constant term $1$. We now derive Good's formula by our
residue theorem:

Let $K\ll x_1,x_2,\ldots ,x_n\gg $ be the working field. Then
$x_i$ is the initial term of $F_i$, and the Jacobian number
$j(F_1,\dots ,F_n)$ equals $1$. Let $y_i=F_i(\mathbf{x})$. We will
have $x_i=G_i(\mathbf{y})$. Then
\begin{align}\label{e-derive-good}
[y_1^{k_1}\cdots y_n^{k_n}] G_i(\mathbf{y}) &= \res_\mathbf{y} y_1^{-1-k_1}\cdots y_n^{-1-k_n} G_i(y)\\
&= \res_\mathbf{x} F_1^{-1-k_1}\cdots F_n^{-1-k_n} x_i
J(\mathbf{F}).
\end{align}
The above argument works the same way by using Jacobi's residue
formula.

 A similar
computation applies to the non-diagonal case by working in $K^\rho
\ll x_1,\dots, x_n\gg$, where $\rho$ is the injective homomorphism
into $K\ll x_1,\dots,x_n,t\gg$ induced by $\rho: x_i\mapsto x_i
t$. This total ordering makes $x_i$ the initial term of $F_i$ for
all $i$, and clearly $K^\rho \ll x_1,\dots, x_n\gg$  contains
$K[[x_1,\dots, x_n]]$ as a subring. This way is equivalent to the
homogeneous expansion introduced in \citep{reversion}. Note that
Jacobi's formula does not apply directly, though Gessel
\cite{gessel-res} showed how the non-diagonal case could be derived
from the diagonal case. Note also that we cannot apply the residue
theorem in $K\ll x_1,\ldots, x_n\gg $, because the Jacobian number
might equal $0$. {}For example, if $x_n$ does not divide $F_n$, then
it is easily seen that the exponent of $x_n$ in the initial term
of $F_i$ is zero for all $i$. So the Jacobian number of $F_1,\dots
,F_n$ is $0$.

More generally, let $\Phi \in K[[y_1,\ldots ,y_n]]$. Then
\begin{equation*}
[y_1^{k_1}\cdots y_n^{k_n}] \Phi(\mathbf{G}(\mb{y}))
=\res_\mathbf{x} F_1^{-1-k_1}\cdots F_n^{-1-k_n} \Phi(\mathbf{x})
J(\mathbf{F}).
\end{equation*}

Multiplying both sides of the above equation by $y_1^{k_1}\cdots
y_n^{k_n}$, and summing on all nonnegative integers
$k_1,k_2,\ldots ,k_n$, we get
\begin{align}\label{e-3-another}
\Phi(\mathbf{G}(\mathbf{y}))=\res_\mathbf{x}
\frac{1}{F_1-y_1}\cdots \frac{1}{F_n-y_n}J(\mathbf{F})
\Phi(\mathbf{x}),
\end{align}
which is true as power series in the $y_i$'s.

It's natural to ask if we can get this formula directly from our
residue theorem. The answer is yes. The argument is given as
follows.

The working field is $K^{\rho}\ll x_1,\ldots ,x_n\gg\ll y_1,\ldots
,y_n\gg $. We let $z_i=F_i-y_i$. Then $x_i=G_i(\mb{y}+\mb{z})$,
and the initial term of $F_i-y_i$ is $x_i$, for $y_i$ has higher
order than the $x$'s. Thus the Jacobian number is $1$, and the
Jacobian determinant $J(\mb{z}|\mb{x})$ still equals to
$J(\mathbf{F})$. Applying the residue theorem, we get
$$\res_\mathbf{x} \frac{1}{F_1-y_1}\cdots \frac{1}{F_n-y_n}J(\mathbf{F})
\Phi(\mathbf{x}) =\res_\mathbf{z} \frac{1}{z_1z_2\cdots z_n}
\Phi(\mathbf{G}(\mb{y}+\mb{z})).$$ Since
$\Phi(\mathbf{G}(\mb{y}+\mb{z}))$ is in $K[[\mb{y},\mb{z}]]$, the
final result is obtained by setting $\mb{z}=\mb{0}$ in
$\Phi(\mathbf{G}(\mb{y}+\mb{z}))$.

Note that $J(\mathbf{F})\in K[[\mb{x}]]$ has constant term $1$.
Therefore $J(\mathbf{F})^{-1}\Phi(\mb{x})$ is also in
$K[[\mb{x}]]$. Hence we can reformulate \eqref{e-3-another} as
$$\res_\mb{x} \frac{1}{F_1-y_1}\cdots \frac{1}{F_n-y_n}\Phi(\mb{x})=
\Phi(\mb{x})J(\mathbf{F})^{-1}|_{\mb{x}=\mb{G}(\mb{y})}.$$

\section{Dyson's Conjecture\label{s-dyson}}
Our residue theorem can be used to prove a conjecture of Dyson.
\begin{thm}[Dyson's Conjecture]\label{t-dyson}
Let $a_1,\ldots ,a_n$ be $n$ nonnegative integers. Then the
following equation holds as Laurent polynomials in $\mb{z}$.
\begin{equation}
\CT_{\mb{z}} \prod_{1\le i\ne j \le n}
\left(1-\frac{z_i}{z_j}\right)^{\!\!a_j} =
 \frac{(a_1+a_2+\cdots a_n)!}{a_1!\, a_2!\, \cdots a_n!}.\label{e-dyson-new}
\end{equation}
\end{thm}

{}For $n=3$ this assertion is equivalent to the familiar Dixon
identity:
\begin{equation*}
\sum_{j}(-1)^j \binom{a+b}{a+j}\binom{b+c}{b+j}\binom{c+a}{c+j}
=\frac{(a+b+c)!}{a!\, b!\, c!}.
\end{equation*}
Theorem \ref{t-dyson} was first proved by Wilson \cite{wilson} and
Gunson \cite{gunson} independently. A similar proof was given by
Egorychev in \cite[p.~151--153]{ego}. These proofs use integrals
of analytic functions. A simple induction proof was found by Good
\cite{good1}. We are going to give a Laurent series proof by using
the residue theorem for MN-series. Our new proof uses Egorychev's
change of variables, and uses Wilson's argument for evaluating the
log Jacobian. This leads to a generalization of Theorem
\ref{t-dyson}.

Let $\mb{z}$ be the vector $(z_1,z_2,\ldots, z_n)$. If $\mb{z}$
appears in the computation, we use $\mb{z}$ for the product
$\mb{z^1}=z_1z_2\cdots z_n$. We use similar notation for $\mb{u}$.

Let $\Delta(\mb{z})=\Delta (z_1,\ldots
,z_n)=\prod_{i<j}(z_i-z_j)=\det(z_i^{n-j})$ be the Vandermonde
determinant in $\mb{z}$, and let
$\Delta_j(\mb{z})=\Delta(z_1,\ldots ,\hat{z}_j,\ldots ,z_n)$,
where $\hat{z}_j$ means to omit $z_j$. We introduce new variables
$u_j=(-1)^{j-1}z_j^{n-1}\Delta_j(\mb{z})$. Then they satisfy the
equations
$$\Delta(\mb{z})=\sum_{j=1}^n (-1)^{j-1} z_j^{n-1}\Delta_j(\mb{z})=u_1+u_2+\cdots +u_n,$$
and
$$u_1\cdots u_n=\prod_{j=1}^n  (-1)^{j-1}z_j^{n-1}\Delta_j(\mb{z}) =(-1)^{\binom{n}{2}}
\mb{z}^{n-1} (\Delta(\mb{z}))^{n-2}.$$ We also have
$$\prod_{i=1,i\ne j}^n \left(1-\frac{z_i}{z_j}\right) =(-1)^{j-1}
\frac{\Delta(\mb{z})}{z_j^{n-1}\Delta_j(\mb{z})}=\frac{u_1+u_2+\cdots
+u_n}{u_j}.$$

Thus equation \eqref{e-dyson-new} is equivalent to
$$\CT_{\mb{z}} \frac{(u_1+u_2+\cdots +u_n)^{a_1+a_2+\cdots +a_n}}{u_1^{a_1}\cdots u_n^{a_n}}
=\frac{(a_1+a_2+\cdots a_n)!}{a_1!\, a_2!\, \cdots a_n!},$$ which
is a direct consequence of the multinomial theorem and the
following proposition.
\begin{prop}\label{p-dyson1}
{}For any series $\Phi(\mb{z})\in K^\mb{u}\ll \mb{z} \gg$, we have
$$\ct_{\mb{z}} \Phi(u_1,\dots ,u_n) =\ct_{\mb{u}} \Phi(u_1,\dots ,u_n).$$
\end{prop}

%

In fact, we can prove a more general formula. Let $r$ be an
integer and let $u^{(r)}_j = (-1)^{j-1}z_j^{r}\Delta_j(\mb{z})$.
Then $u^{(r)}_1+\cdots +u^{(r)}_n$ equals $ h_{r-n+1}
(z_1,z_2,\dots ,z_n)\Delta(\mb{z})$ for $r\ge n-1$ and equals $0$
for $0\le r\le n-2$, where $h_k(\mb{z})=\sum_{i_1\le \cdots \le
i_k}z_{i_1}\cdots z_{i_k}$ is the complete symmetric function
\cite[Theorem 7.15.1]{EC2}. We have the following generalization.

\begin{thm}\label{t-dyson-g}
If $r$ is not equal to any of $0,1,\cdots ,n-2,$ or
$-\binom{n-1}{2}$, then for any series $\Phi(\mb{z})\in K^\rho\ll
\mb{z} \gg$, where $\rho(z_i)= u_i^{(r)}$, we have
$$\ct_{\mb{z}} \Phi(u^{(r)}_1,\dots ,u^{(r)}_n) =\ct_{\mb{u}^{(r)}} \Phi(u^{(r)}_1,\dots
,u^{(r)}_n).$$
\end{thm}
Note that Proposition \ref{p-dyson1} is the special case for
$r=n-1$ of Theorem \ref{t-dyson-g}. If we set $r=n$, the
multinomial theorem yields the following:

\begin{cor}
Let $a_1,\ldots ,a_n$ be $n$ nonnegative integers. Then the
following equation holds for Laurent polynomials in $\mb{z}$:
\begin{equation}
\CT_{\mb{z}} \frac{(z_1+\cdots +z_n)^{a_1+\cdots
+a_n}}{z_1^{a_1}\cdots z_n^{a_n}}\prod_{1\le i\ne j \le n}
\left(1-\frac{z_i}{z_j}\right)^{\!\!a_j} =
 \frac{(a_1+a_2+\cdots +a_n)!}{a_1!\, a_2!\, \cdots a_n!}.\label{e-dyson}
\end{equation}
\end{cor}

By Theorem \ref{t-MNresidue} and Proposition \ref{p-3-ljacobian},
Theorem \ref{t-dyson-g} is equivalent to the assertion that the
log Jacobian is a nonzero constant. To show this, we use

\begin{lem}[\cite{wilson}, Lemma 4]\label{l-dyson-sym}
Let $G(x_1,\dots ,x_n)$ be a ratio of two polynomials in the
$x$'s, in which the denominator is $\Delta(x_1,\dots ,x_n)$
{and}
\begin{enumerate}
\item $G$ is a symmetric function of $x_1,\dots ,x_n$.
\item $G$ is homogeneous of degree $0$ in the $x$'s.
\end{enumerate}
Then $G$ is a constant.
\end{lem}

\begin{proof}[Proof of Theorem \ref{t-dyson-g}]
In order to compute the log Jacobian, we let
$$J=\det(J_{ij})=\det \left( \frac{\partial \log u^{(r)}_i}{\partial \log z_j}\right).$$ Then
$J_{ii}=r$ and $J_{ij}=\sum_{k\ne i} \frac{z_i}{z_k-z_j}$ for
$i\ne j$. We first show that $J$ is a constant by Lemma
\ref{l-dyson-sym}. It is easy to see that $J$ satisfies conditions
$1,2$ in Lemma \ref{l-dyson-sym}. Now we show that the denominator
of $J$ is $\Delta(\mb{z})$, so that we can claim that the Jacobian
is a constant, and hence equals the Jacobian number.

Evidently $J$ is the ratio of two polynomials in the $\mb{z}$'s,
whose denominator is a product of factors $z_i-z_j$ for some $i\ne
j$. {}From the expression of $J_{ij}$, we see that $z_i-z_j$ only
appears in the $i$th and the $j$th column. Every $2$ by $2$ minor
of the $i$th and $j$th  columns is of the following form, in which
we assume that $k$ and $l$ are not one of $i$ and $j$.
$$\left|\begin{array}{cc}
J_{ki} & J_{kj} \\
J_{li} & J_{lj}
\end{array} \right|
=\left|\begin{array}{cc} \frac{z_k}{z_j-z_i}+\sum_{s\ne i,j}
\frac{z_k}{z_s-z_i} & \frac{z_k}{z_i-z_j}+\sum_{s\ne i, j}
\frac{z_k}{z_s-z_j} \\
\frac{z_l}{z_j-z_i}+ \sum_{s\ne i,j} \frac{z_l}{z_s-z_i} &
\frac{z_l}{z_i-z_j} +\sum_{s\ne i, j} \frac{z_l}{z_s-z_j}
\end{array} \right|.
$$
In the above determinant, the terms containing $(z_i-z_j)^2$ as
the denominator cancel. Therefore, expanding the determinant
according to the $i$th and  $j$th column, we see that
$\Delta(\mb{z})$ is the denominator of $J$.

Now the initial term of $z_i-z_j$ is $z_i$ if $i<j$. We see that
the initial term of $u^{(r)}_1$ is $z_1^rz_2^{n-2}z_3^{n-3}\cdots
z_{n-1}$. Similarly we can get the initial term for $u^{(r)}_j$.
The Jacobian number, denoted by $j(r)$, is thus the determinant
$$j(r)=\det \left(
\begin{array}{ccccc}
r & n-2 & n-3&  \cdots &0 \\
n-2 &r &  n-3&  \cdots & 0 \\
\vdots &\vdots &\vdots  & \vdots &\vdots \\
n-2 & n-3 & n-4 & \cdots & r
\end{array}
\right),
$$
where the displayed matrix has diagonal entries $r$, and other
entries in each row are $n-2,n-3,\dots , 0$, from left to right.

Since the row sum of each row is $r+\binom{n-1}{2}$, it follows
that $j(-\binom{n-1}{2})=0$. We claim that $j(r)=0$ when
$r=0,1,\dots , n-2$. {}For in those cases, $u_1^{(r)}+\cdots
+u_{n}^{(r)} =0$. This implies that the Jacobian is $0$, and hence
$j(r)=0$. We can regard $j(r)$ as a polynomial in $r$ of degree
$n$, and we already have $n$ zeros. So $j(r)=r(r-1)\cdots (r-n+2)
(r+\binom{n-1}{2})$ up to a constant. This constant equals $1$
through comparing the leading coefficient of $r$.

In particular, $j(n-1)= \binom{n}{2} (n-1)!= (n-1) n!/2.$ Note
that in \cite[p.~153]{ego}, the constant was said to be
$(2n-3)(n-1)!$, which is not correct.
\end{proof}

Another proof of Dyson's conjecture by our residue theorem is to
use the change of variables by Wilson \citep{wilson}.

Let
$$v_j=\prod_{ i=1, i\ne j}^n (1-z_j/z_i)^{-1}.$$
Then the initial term of $v_j$ is $z_j^{-(n-j)}z_{j+1}\cdots z_n$
up to a constant. Since the order of $v_n$ is $\mb{0}$, we have to
exclude $v_n$ from the change of variables, for otherwise, the
Jacobian number will be $0$. In fact, we have the relation
$v_1+v_2+\cdots +v_n=1$, which can be easily shown by Lemma
\ref{l-dyson-sym}.

Dyson's conjecture is equivalent to
\begin{align}\label{e-3-dyson-v}
\ct_{\mb{z}} \prod_{j=1}^n v_i^{-a_j} = \frac{(a_1+a_2+\cdots
a_n)!}{a_1!a_2!\cdots a_n!}
\end{align}

\begin{proof}[Another Proof of Dyson's Conjecture]
Using Lemma \ref{l-dyson-sym} and Wilson's argument, we can
evaluate the following log Jacobian. (See \citep{wilson} for
details.)
$$\frac{\partial (\log v_1,\log v_2 ,\dots ,\log v_{n-1})}{\partial(
\log z_1,\log z_2,\dots ,\log z_{n-1})}=(n-1)!\, v_n.$$

Then by the residue theorem
$$\ct_z \Phi(v_1,\dots ,v_{n-1},z_n) =\ct_{v_1,\dots,v_{n-1},z_n}
(1-v_1-\cdots -v_{n-1})^{-1} \Phi(v_1,\dots, v_{n-1},z_{n}).$$

In particular (since the initial term of $1-v_1-\cdots -v_{n-1}$
is $1$) we have:
\begin{align*}
\ct_{\mb{z}} \prod_{j=1}^n v_i^{-a_j}
&=\ct_{v_1,\dots,v_{n-1},z_n}
(1-v_1-\cdots -v_{n-1})^{-a_n-1}\prod_{j=1}^{n-1} v_i^{-a_j}\\
&=[v_1^{a_1}\cdots v_{n-1}^{a_{n-1}}] \sum_{m\ge 0 }
\binom{a_n+m}{a_n}
(v_1+\dots +v_{n-1})^m\\
&=\binom{a_n+a_1+\cdots + a_{n-1}}{a_n} \binom{a_1+\cdots +
a_{n-1}}{a_1,\dots,a_{n-1}}.
\end{align*}
Equation \eqref{e-3-dyson-v} then follows.
\end{proof}

\vspace{3mm}\noindent {\bf Acknowledgement:} The author is very
grateful to Ira Gessel's guidance and help.

\bibliographystyle{amsplain}

\providecommand{\bysame}{\leavevmode\hbox
to3em{\hrulefill}\thinspace}
\providecommand{\MR}{\relax\ifhmode\unskip\space\fi MR }
\providecommand{\MRhref}[2]{%
  \href{http://www.ams.org/mathscinet-getitem?mr=#1}{#2}
} \providecommand{\href}[2]{#2}

\end{document}